# Smoothing-inspired lack-of-fit tests based on ranks

Jeffrey D. Hart[*,1]

*Texas A&M University*

**Abstract:** A rank-based test of the null hypothesis that a regressor has no effect on a response variable is proposed and analyzed. This test is identical in structure to the order selection test but with the raw data replaced by ranks. The test is nonparametric in that it is consistent against virtually any smooth alternative, and is completely distribution free for all sample sizes. The asymptotic distribution of the rank-based order selection statistic is obtained and seen to be the same as that of its raw data counterpart. Exact small sample critical values of the test statistic are provided as well. It is shown that the Pitman-Noether efficiency of the proposed rank test compares very favorably with that of the order selection test. In fact, their asymptotic relative efficiency is identical to that of the Wilcoxon signed rank and $t$-tests. An example involving microarray data illustrates the usefulness of the rank test in practice.

## 1. Introduction

It is arguable that the best and most aesthetically appealing ideas in science are those that combine the virtues of simplicity and effectiveness. Perhaps no idea in statistics better achieves this combination than rank-based methods. In this paper such methods are brought to bear on the problem of nonparametrically testing lack-of-fit in regression.

Nonparametric lack-of-fit tests based on smoothing methods have received a great deal of attention in recent years. Many of these tests are discussed by Hart [3]. An omnipresent problem of smoothing-based tests, as with any other test, is uncertainty about the sampling distribution of the test statistic. The bootstrap seems to have become the principal way of dealing with this problem when the distribution of the data is unknown. However, the bootstrap is not a panacea for at least two reasons. First of all, its performance usually breaks down when the data do not possess at least two moments, and secondly, even if the moment assumptions are met, it only guarantees test validity asymptotically, as the number of data increase without bound.

Rank tests have long been an effective way of dealing with uncertainty about the data distribution. An important, but not nearly exhaustive, set of references on the subject is Chernoff and Savage [1], Puri and Sen [14], Randles and Wolfe

---

[*]Supported by NSF Grants DMS-02-03884 and DMS-06-04801.
[1]Department of Statistics, Texas A&M University, College Station, TX 77843-3143, USA, e-mail: hart@stat.tamu.edu

[15], Hettmansperger and McKean [5]. Rank tests are relatively insensitive to assumptions about the underlying data distribution. Indeed, under quite general assumptions they are often completely distribution-free, no matter the sample size. Furthermore, they usually sacrifice little, if any, power relative to analogous tests based on the raw data. The purpose of this paper is to consider certain smoothing-based lack-of-fit tests when they are applied to ranks of the responses rather than the responses themselves. It will be shown that many of the conclusions about rank tests from other settings carry over to the nonparametric lack-of-fit problem.

Aside from Hart [3], it appears that the work of Lombard [9] is more closely related to that in the current paper than any other in the existing literature. Lombard [9] considers the problem of testing a sequence of independent random variables for mean or scale constancy. He proposes that the periodogram of the ranks of the data be computed, and defines a test statistic as a weighted sum of centered periodogram ordinates. The omnibus nature of this statistic and the fact that it uses Fourier transforms of ranks make it similar in spirit to statistics proposed in the current paper.

Many, and probably most, areas of applied research require good methods for testing the fit of parametric regression models. The methodology proposed in this paper thus has the potential of being extremely useful in interdisciplinary research. Evidence of this potential is provided in Section 4, which involves an example from the field of microarray analysis. Our methods are used to establish that test samples differ from reference samples in terms of DNA copy number. Such differences are known to be correlated with cancer incidence in the test subjects.

In the next section we will introduce a nonparametric regression model and propose some rank-based tests of the null hypothesis that the regression function is constant. Asymptotic properties of one of these tests are developed in this section as well. In Section 3 we discuss how the basic ideas from Section 2 may be extended to test the fit of a linear model. The aforementioned illustration involving microarrays is the subject of Section 4, concluding remarks are made in Section 5, and proofs of three theorems are given in an Appendix.

## 2. Testing the no-effect hypothesis

Suppose that one observes responses $Y_1, \ldots, Y_n$ at fixed design points $x_1, \ldots, x_n$ which, mainly for convenience, we take to be $x_i = (i - 1/2)/n$, $i = 1, \ldots, n$. In general it is assumed that

(2.1) $$Y_i = r(x_i) + \epsilon_i, \quad i = 1, \ldots, n,$$

where $r$ is some function that is square integrable over $[0, 1]$ and $\epsilon_1, \ldots, \epsilon_n$ are independent and identically distributed. A fundamentally important problem in this setting is establishing that the response is indeed related to the design variable. Formally, we wish to test the null hypothesis

$$H_0 : r \equiv C,$$

against the alternative

$$H_1 : \int_0^1 (r(x) - \bar{r})^2 \, dx > 0,$$

where $C$ is an unknown constant and $\bar{r} = \int_0^1 r(x) \, dx$.

A nonparametric means of estimating $r$ is to use an orthogonal series $\hat{r}_m$ of the form

$$(2.2) \qquad \hat{r}_m(x) = \begin{cases} \tilde{\phi}_0, & m = 0 \\ \tilde{\phi}_0 + 2\sum_{j=1}^{m} \tilde{\phi}_j \cos(\pi j x), & 1 \leq m < n, \end{cases}$$

where

$$\tilde{\phi}_j = \frac{1}{n}\sum_{i=1}^{n} Y_i \cos(\pi j x_i), \quad j = 0, \ldots, n-1.$$

The integer $m$ is the smoothing parameter of $\hat{r}_m$. Many tests of $H_0$ vs. $H_1$ based on Fourier series smoothers have been proposed; see Chapter 7 of Hart [3] for a description of a number of these. One such test is the so-called *order selection* (OS) test of Eubank and Hart [2]. Define

$$T_n = \max_{1 \leq m < n} \frac{1}{m} \sum_{j=1}^{m} \frac{2n\tilde{\phi}_j^2}{\hat{\sigma}^2},$$

where $\hat{\sigma}^2$ is a consistent estimator of $\sigma^2 = \text{Var}(\epsilon_1)$, assuming that it exists finite. The OS test rejects $H_0$ for large values of $T_n$. This test derives its name from the fact that rejecting $H_0$ if and only if $T_n \geq A$ is equivalent to rejecting $H_0$ if and only if the maximizer of

$$(2.3) \qquad M_n(m) = \begin{cases} 0, & m = 0 \\ \sum_{j=1}^{m} 2n\tilde{\phi}_j^2/\hat{\sigma}^2 - Am, & m = 1, 2, \ldots, n-1 \end{cases}$$

with respect to $m$ is greater than 0. When $A = 2$, $M_n$ is precisely the well-known Mallows' criterion (Mallows [11]) for choosing the order of the smoother $\hat{r}_m$.

A rank-based analog of any series-based test is easily defined by using ranks instead of $Y_1, \ldots, Y_n$ in the coefficients $\tilde{\phi}_1, \ldots, \tilde{\phi}_{n-1}$. Let $R(Y_i)$ denote the rank of $Y_i$ among $Y_1, \ldots, Y_n$, and define

$$U_i = \frac{R(Y_i)}{n+1}, \quad i = 1, \ldots, n,$$

and

$$(2.4) \qquad \hat{\phi}_j = \frac{1}{n}\sum_{i=1}^{n} U_i \cos(\pi j x_i), \quad j = 1, \ldots, n-1.$$

A rank-based analog of $T_n$ is

$$R_n = \max_{1 \leq m < n} \frac{1}{m} \sum_{j=1}^{m} \frac{2n\hat{\phi}_j^2}{1/12},$$

which was first proposed in Hart [3]. Here we shall verify properties of $R_n$ that were only conjectured in Hart [3].

### 2.1. Null distribution of the test statistic

When the null hypothesis is true, the distribution of $R_n$ is completely independent of the distribution of $\epsilon_1$. For small values of $n$ one may determine the distribution of $R_n$ exactly by tabulating the value of $R_n$ for each permutation of the ranks $1, \ldots, n$. The following theorem provides the *asymptotic* distribution of $R_n$.

TABLE 1
*Tail probabilities for $R_n$ at large sample quantiles. Values for $5 \leq n \leq 10$ are exact, while those for $n = 15$, 20 and 30 are gotten from simulation. The latter values are accurate to within 0.0005 with 95% confidence*

| | Large sample quantile | | | |
|---|---|---|---|---|
| $n$ | **3.221** | **4.179** | **6.745** | **10.850** |
| 5 | 0.1000 | 0.0167 | 0.0000 | 0.0000 |
| 6 | 0.1028 | 0.0417 | 0.0000 | 0.0000 |
| 7 | 0.1040 | 0.0476 | 0.0004 | 0.0000 |
| 8 | 0.1034 | 0.0487 | 0.0022 | 0.0000 |
| 9 | 0.1042 | 0.0482 | 0.0039 | 0.0000 |
| 10 | 0.1030 | 0.0485 | 0.0053 | 0.0000 |
| 15 | 0.1030 | 0.0496 | 0.0078 | 0.0002 |
| 20 | 0.1020 | 0.0496 | 0.0086 | 0.0003 |
| 30 | 0.1016 | 0.0501 | 0.0089 | 0.0006 |
| $\infty$ | 0.10 | 0.05 | 0.01 | 0.001 |

**Theorem 2.1.** *Suppose that model (2.1) holds with $r \equiv$ constant and $\epsilon_1, \ldots, \epsilon_n$ independent and identically distributed as $F$, where $F$ is a continuous cumulative distribution function (cdf). Then*

$$(2.5) \qquad \lim_{n \to \infty} P(R_n \leq t) = \exp\left\{-\sum_{j=1}^{\infty} \frac{P(\chi_j^2 > jt)}{t}\right\} \equiv G(t),$$

*where $\chi_j^2$ has the chi-squared distribution with $j$ degrees of freedom, $j = 1, 2, \ldots$.*

The distribution $G$ is precisely the same limiting distribution that $T_n$ has under appropriate moment conditions (Chapter 7, Hart [3]). The minimum number of moments that $\epsilon_1$ must possess in order for $T_n$ to have limiting distribution $G$ is two. This in itself provides a key motivation for use of the rank-based test, since (2.5) requires no moment conditions whatsoever.

Table 1 gives tail probabilities for $R_n$ at selected large sample percentiles. It seems clear from this table that use of the large sample percentiles will result in a conservative test so long as $\alpha \leq 0.05$. At $\alpha = 0.10$, use of the large sample quantiles results in tests that are only slightly liberal.

### 2.2. Power of rank-based order selection test

Another attractive feature of many rank tests is that they give up remarkably little power in comparison to tests based on the raw data. They can even be more powerful than raw-data counterparts when the parent distribution is sufficiently heavy-tailed. The reader is referred to Puri and Sen [14] for a review of some of these results.

Before investigating power of the rank-based OS test, we can gain some intuition by considering the function

$$\mu(x) = 2\sum_{i=1}^{\infty} \phi_j \cos(\pi j x),$$

where

$$\phi_j = \lim_{n \to \infty} E\hat{\phi}_j, \quad j = 1, 2, \ldots.$$

It should not be surprising that the rank-based OS test has power comparable to that of an ordinary OS test applied to model (2.1) with $r \equiv \mu$ and $\text{Var}(\epsilon_1) = 1/12$. Letting $I(A)$ denote the indicator of event $A$, we have $R(Y_i) = \sum_{j=1}^{n} I(Y_i \geq Y_j)$ and hence

$$
\begin{aligned}
E\hat{\phi}_j &= \frac{1}{n^2} \sum_{i=1}^{n} \sum_{k=1}^{n} P(Y_i \geq Y_k) \cos(\pi j x_i) \\
&= \frac{1}{n^2} \sum_{i=1}^{n} \sum_{k \neq i} H\left(r(x_i) - r(x_k)\right) \cos(\pi j x_i),
\end{aligned}
$$

where $H$ is the cdf of $\epsilon_1 - \epsilon_2$. It follows that

$$
E\hat{\phi}_j = \int_0^1 \mu(x) \cos(\pi j x) \, dx + O(n^{-1}),
$$

where the $O$ term is bounded uniformly in $j$ and

$$
\mu(x) = \int_0^1 H(r(x) - r(u)) \, du.
$$

When the alternative $r$ is very close to the null, the form of $\mu$ relative to $r$ is more transparent. If $r(x) = \beta(x)/\sqrt{n}$, then

$$
E\hat{\phi}_j = \frac{h(0)}{\sqrt{n}} \int_0^1 \beta(x) \cos(\pi j x) \, dx + o(n^{-1/2}),
$$

where $h$ is the density correpsonding to $H$. For the local alternative $\beta/\sqrt{n}$, one thus anticipates that the rank-based OS test will have power comparable to the ordinary OS test against the local alternative $h(0)\beta/\sqrt{n}$ with $\sigma = \sqrt{1/12}$. This indeed turns out to be the case, as we now show.

The relative power of the rank-based and ordinary OS tests will be investigated using the notion of Pitman-Noether efficiency. We consider a sequence of local alternatives to $H_0$ having the form

$$
(2.6) \qquad r_n(x) = \frac{\beta(x)}{\sqrt{n}},
$$

where $\int_0^1 \beta^2(x) \, dx > 0$. Let $P_n(r, f)$ and $P_n^R(r, f)$ denote the power of the ordinary and rank-based OS tests, respectively, when the regression function is $r$, the error density is $f$, the sample size is $n$ and the level of significance is the same for each test. Under appropriate conditions, each of $P_n(r_n, f)$ and $P_n^R(r_n, f)$ has a limit greater than 0 and less than 1 as $n \to \infty$. If we can find a function of $n$, call it $n^*$, such that

$$
\lim_{n \to \infty} P_n(r_n, f) = \lim_{n \to \infty} P_{n^*}^R(r_n, f),
$$

then the Pitman-Noether asymptotic relative efficiency (ARE) of the rank-based OS test to the ordinary OS test is

$$
e(f) = \lim_{n \to \infty} \frac{n}{n^*}.
$$

This notion of efficiency relates how many more (or fewer) observations are needed for a rank-based order selection test to have the same power as the analogous test based on the raw data.

We now state a theorem concerning the limiting power of the rank-based OS test.

**Theorem 2.2.** *Suppose that model (2.1) holds with $r(x) \equiv r_n(x)$, as defined by (2.6). Assume that the following two conditions hold:*

- *The density $h$ (of $\epsilon_1 - \epsilon_2$) is Lipschitz continuous in a neighborhood of 0.*
- *The function $\beta$ is Lipschitz continuous on $[0, 1]$, except perhaps at finitely many points where it is simply discontinuous.*

Let $Z_1, Z_2, \ldots$ be i.i.d. standard normal random variables, $h$ be the density of $\epsilon_1 - \epsilon_2$, and $\beta_j = \int_0^1 \beta(x) \cos(\pi j x)\, dx$, $j = 1, 2, \ldots$. Then, for a sequence of level $\alpha$ tests,

$$\lim_{n \to \infty} P_n^R(r_n, f) = P(\mathcal{T} \geq t_\alpha),$$

where $t_\alpha$ is the $1 - \alpha$ quantile of $G$, the limiting null distribution of $R_n$, and

$$\mathcal{T} = \max_{m \geq 1} \frac{1}{m} \sum_{j=1}^m \left( Z_j + \sqrt{24} h(0) \beta_j \right)^2.$$

Under the conditions in Theorem 2.2 and the additional condition that $\epsilon_1$ has four moments finite, Theorem 7.10, p. 201 of Hart [3] states that

$$(2.7) \quad \lim_{n \to \infty} P_n(r_n, f) = P \left[ \max_{m \geq 1} \frac{1}{m} \sum_{j=1}^m \left( Z_j + \frac{\sqrt{2} \beta_j}{\sigma} \right)^2 \geq t_\alpha \right],$$

where $\sigma^2 = \mathrm{Var}(\epsilon_1)$.

Combining Theorem 2.2 with (2.7) allows us to obtain the ARE of the rank-based to the ordinary OS test.

**Theorem 2.3.** *Let the conditions of Theorem 2.2 hold and suppose that $\epsilon_1$ has four moments finite. Then the ARE of the rank-based OS test to the ordinary OS test is*

$$(2.8) \quad \eta(f) = 12\sigma^2 \left( \int_{-\infty}^{\infty} f^2(x)\, dx \right)^2.$$

Interestingly, the ARE in Theorem 2.3 is precisely the same as those of the Wilcoxon signed-rank test to the $t$-test and the Mann-Whitney test to the two-sample $t$-test in classical versions of the one- and two-sample location problems, respectively (Randles and Wolfe [15]). This fact parallels results of Hettmansperger and McKean [5] in the problem of testing the fit of a linear model. It is shown in Hettmansperger and McKean [5] (pp. 176–178) that the asymptotic efficiencies of rank-based tests relative to the classical $F$-test are the same as corresponding AREs in the one- and two-sample location problems. In both Hettmansperger and McKean [5] and our setting, the ARE turns out to be simply the ratio of noncentrality parameters. In our problem, for example,

$$ARE = \frac{24 h^2(0) \beta_j^2}{2 \beta_j^2 / \sigma^2} = 12 \sigma^2 h^2(0).$$

For any scale family $f(x) = f_0(x/\sigma)/\sigma$ such that $\int_{-\infty}^{\infty} x^2 f_0(x)\, dx = 1$, we have

$$\eta(f) = 12 \left( \int_{-\infty}^{\infty} f_0^2(x)\, dx \right)^2.$$

A special case of interest is when $f_0$ is the standard normal density, $\phi$. It is well-known that $e(\phi) = 3/\pi \approx 0.955$, and hence the loss in power using the rank-based test is quite small.

For long-tailed distributions, the rank test can be more efficient than the ordinary OS test. Let $t_k$ denote the density of a $t$ distribution with $k$ degrees of freedom, where $k$ is a positive integer. Then it can be verified that $e(t_k) > 1$ for $5 \leq k \leq 18$, with $e(t_5) \approx 1.24$. Of interest are empirical studies to determine how large $n$ must be in order for such power improvements to be realized. Results of McKean and Sheather [12] show that rank-based $F$-tests for testing the fit of linear models are indeed more powerful in small samples than an ordinary $F$-test when the error distribution is heavier tailed than the normal.

Borrowing an idea from the theory of linear rank statistics, one could replace $U_i$ in (2.4) by a normal score, i.e., $\Phi^{-1}(U_i)$, where $\Phi$ is the standard normal cdf. Chernoff and Savage [1] show that using normal scores in the classical two-sample location problem yields an ARE, relative to the $t$-test, that is always at least 1. One anticipates a similar benefit from using normal scores in the setting of the current paper, but we do not pursue this question further here.

### 2.3. Other series-based rank tests

As mentioned previously, a number of lack-of-fit tests based on orthogonal series have been proposed. Rank-based analogs of any of these can be performed simply by replacing the raw data by ranks. Perhaps the most fundamental series-based test is the regression analog of a Neyman smooth test (Neyman [13]). The basis for such a test is the statistic

$$S_{n,m} = \sum_{j=1}^{m} \frac{2n\tilde{\phi}_j^2}{\hat{\sigma}^2},$$

where $m$ is fixed prior to data collection. This statistic has a simple and familiar asymptotic null distribution, i.e., chi-squared with $m$ degrees of freedom. According to Lehmann [8], it also has a uniformly most powerful property when the data are normally distributed and the function $r$ has the form

$$r(x) = a_0 + 2\sum_{j=1}^{m} a_j \cos(\pi j x).$$

A difficulty with $S_{n,m}$ is its dependence on $m$, a poor choice of which could result in a loss of power or even inconsistency. Ledwina [7] and Kuchibhatla and Hart [6] proposed that one use a test statistic of the form $S_{n,m}$ with a data-driven choice for $m$. Here we define a rank-based analog of such a statistic. First define the Mallows-like criterion $M_n^R$:

$$M_n^R(m) = \begin{cases} 0, & m = 0 \\ \sum_{j=1}^{m} 24n\hat{\phi}_j^2 - 2m, & m = 1, \ldots, n-1. \end{cases}$$

If $\hat{m}$ is the maximizer of $M_n^R$, then we may define

$$S_{n,\hat{m}}^R = \sum_{j=1}^{\hat{m}} 24n\hat{\phi}_j^2,$$

and reject $H_0$ for large values of $S_{n,\hat{m}}^R$. This statistic is distribution-free under $H_0$

and its critical values are easily approximated by simulation. One could also use a BIC-like criterion in place of $M_n^R$ by replacing the term $2m$ by $(\log n)m$. The criterion, Mallows or BIC, used to choose $m$ in $S_{n,m}$ has a nontrivial effect on properties of the resulting tests. The reader is referred to Section 7.7.4 of Hart [3] for a discussion of these properties.

Hart [4] proposed a Bayesian-motivated lack-of-fit test with very good overall power properties. A rank-based analog of this statistic is

$$B_n = \sum_{j=1}^{n-1} j^{-2} \exp\left(12n\hat{\phi}_j^2\right).$$

Under $H_0$, $B_n$ is distribution-free and converges in distribution to the random variable $\sum_{j=1}^{\infty} j^{-2} \exp(Z_j^2/2)$ as $n \to \infty$, where $Z_1, Z_2, \ldots$ are i.i.d. standard normal. Under appropriate regularity conditions the ARE of a test based on $B_n$ relative to that of the analogous raw-data test (as in Hart [4]) is equal to (2.8). This means that for normal and many other error distributions $B_n$ will inherit all the desirable power properties discussed in Hart [4].

## 3. Testing the fit of linear models

Any of the rank tests previously discussed can be used to test the fit of a linear model. Here one wishes to test a hypothesis of the form

(3.1) $$H_0 : r(x) = \sum_{j=1}^{p} \theta_j r_j(x),$$

where $r_1, \ldots, r_p$ are known functions and $\theta_1, \ldots, \theta_p$ unknown parameters. Rather than ranking the observed data, one ranks residuals from the fitted linear model. Otherwise the tests are done in precisely the same way as before. Now, one is effectively testing the null hypothesis that the expected value of each residual is 0. An excellent reference for rank tests based on residuals in a linear models setting is Hettmansperger and McKean [5].

There is, however, an important difference between testing for constancy of $E(Y_i)$ and testing the fit of a linear model. In the latter problem the joint distribution of the ranks of residuals depends, in general, upon the error distribution. Under the null hypothesis (3.1), the residuals have the form

$$\begin{aligned} e_i &= Y_i - \hat{Y}_i \\ &= \epsilon_i + \sum_{j=1}^{p}(\theta_j - \hat{\theta}_j)r_j(x_i), \quad i = 1, \ldots, n. \end{aligned}$$

If $\hat{\theta}_j$ is a $\sqrt{n}$-consistent estimator of $\theta_j$ for each $j$ (as, for example, least squares estimators generally would be), then the asymptotic distribution of the OS rank test applied to residuals will generally have the same limit distribution as before, i.e., (2.5). However, in small samples the distribution of that statistic can depend upon the error distribution, due to the terms $\sum_{j=1}^{p}(\theta_j - \hat{\theta}_j)r_j(x_i)$, $i = 1, \ldots, n$.

A main competitor of rank methods for dealing with uncertainty about the error distribution is the bootstrap. In the linear models context each of these methods yields valid tests only in the limit, as $n \to \infty$. An interesting question is which of the two methods tends to require larger sample sizes for validity? Intuitively it seems that the rank test would be less sensitive to the underlying error distribution than is the bootstrap. After all, in the absence of any error in estimating the regression coefficients, the rank test would be completely distribution free, whereas a bootstrapped OS test (using the raw data) would not be. Nonetheless, we make no claims in this regard and postpone a comparison of the bootstrap and rank methods to future research.

## 4. A data analysis

Here we apply the rank-based OS test to data collected by the authors of Snijders et al. [16]. The data are from a microarray experiment that measured genome-wide DNA copy number. The variable considered here is the ratio of dye intensities for test and reference samples at a given marker along a chromosome of interest. Each intensity is proportional to the number of marker copies. The reference samples are diploid, and hence each reference marker has only two copies. It is of interest to detect regions on a chromosome where the test samples may have more or fewer copy numbers than the corresponding reference samples. Such variations in copy number are known to be correlated with cancer incidence; see, for example, Lucito et al. [10].

Data sets for four different chromosomes (gotten from cell line GM03563) are shown in Figure 1. In each graph, the horizontal axis is marker location and the vertical axis is the normalized average of three readings of $\log_2(\rho)$, where $\rho$ is the aforementioned ratio of intensities (test over reference). Both the ordinary and rank versions of the OS test were applied to the four data sets to test for constancy of expected $\log_2(\rho)$. In fact, the hypothesis of interest is that this expectation is identically 0. The OS tests have no power against a simple shift alternative, and hence if the null hypothesis of constancy is not rejected, then one would still want to investigate the (perhaps unlikely) possibility that the true function is identical to a nonzero constant.

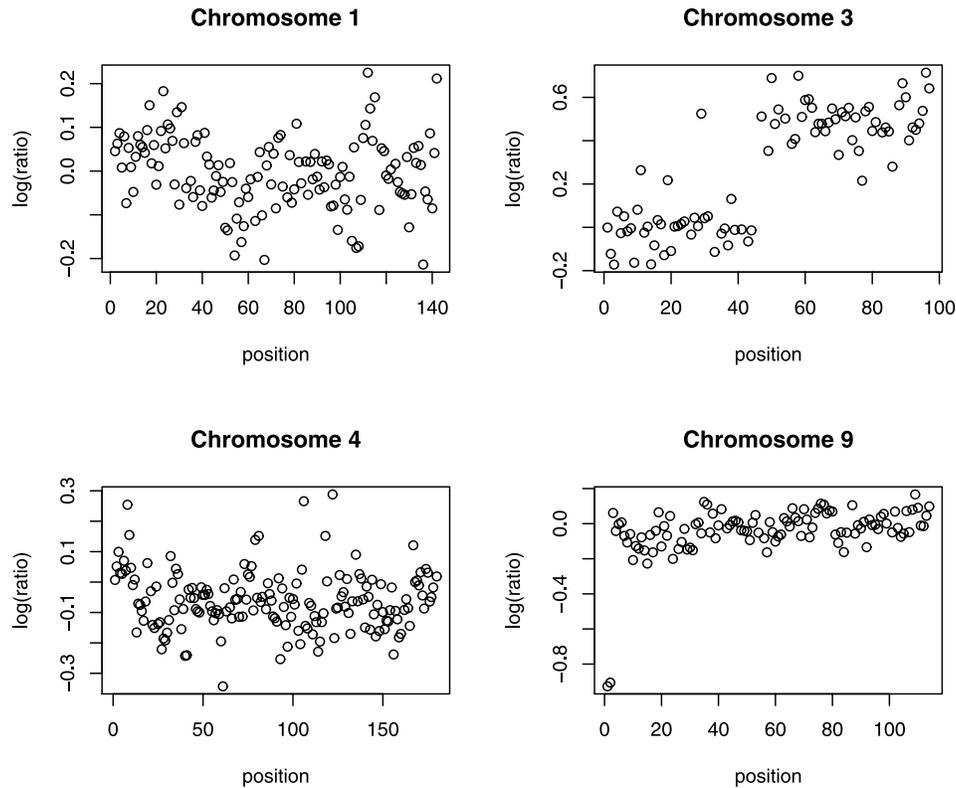

FIG 1. *Microarray data comparing DNA copy number in test and reference samples.*

TABLE 2
*Test statistics and P-values for rank and ordinary OS tests. The values in columns 3 and 5 are large-sample P-values for the rank and ordinary OS tests, respectively*

| Data set | $R_n$ | $P$-value | $T_n$ | $P$-value |
|---|---|---|---|---|
| Chromosome 1 | 8.44 | 0.00379 | 11.05 | 0.00089 |
| Chromosome 3 | 47.50 | 0.00000 | 306.92 | 0.00000 |
| Chromosome 4 | 5.26 | 0.02465 | 5.77 | 0.01795 |
| Chromosome 9 | 17.44 | 0.00003 | 38.99 | 0.00000 |

Values of test statistics and large sample $P$-values are given in Table 2. Since the sample sizes for chromosomes 1, 3, 4 and 9 are 135, 85, 171 and 109, respectively, large sample tests seem more than adequate. Both the ordinary and rank-based OS tests are significant at the 0.05 level of significance for all four data sets, indicating that there are differences between test and reference samples for all four chromosomes. An idea of how test and reference differ can be gotten by considering Fourier series smooths of the form (2.2). It is also interesting to consider the function estimated by a smooth of this same form but with $\tilde{\phi}$ replaced by $\hat{\phi}$. This smooth estimates the function $\mu(x)$, which determines the power of the rank-based test. Estimates of $r(x)/\sigma$ and $\sqrt{12}\mu(x)$ are given in Figure 2. In each of the four plots, the truncation point, $m$, of the estimate $\hat{r}_m$ is chosen by criterion (2.3) with $A = 4.18$, and for simplicity the other truncation point is taken to be the same. (In the case of chromosome 9, two outliers were omitted in computing the smooths.) It is interesting how well the shapes of the rank-based smooths mimic those of the raw-data smooths. For these data, the rank-based smooths have somewhat smaller amplitudes, which is consistent with the rank statistics being smaller than the ordinary OS statistics. Of course, this need not be the case, since, as pointed out earlier, the rank tests will sometimes be more powerful.

It is worthwhile mentioning that the rank-based OS tests are essentially impervious to the two outliers in the chromosome 9 data, evident in Figure 1. It is thus not necessary to delete these cases to determine their effect on the question of whether or not the underlying curve is constant. In contrast, in applying the ordinary OS test to the whole data set, one would wonder whether significance, or lack thereof, was caused by these two points alone. (In this case it turns out that the ordinary OS test is highly significant whether or not the two points in question are included.)

## 5. Conclusions

Nonparametric tests of the null hypothesis that a response and a regressor are unrelated have been proposed and analyzed. The tests are rank-based versions of the order selection test and are completely distribution free under the single condition that the regression errors are independent and identically distributed. These tests have the same surprisingly good power properties possessed by rank tests in simpler testing problems. Rank-based versions of other smoothing-inspired lack of fit statistics were also discussed. One such test promises to have better overall power properties than the rank-based order selection test. The widespread nature of regression analysis suggests that these results could have a substantial impact on a number of different fields.

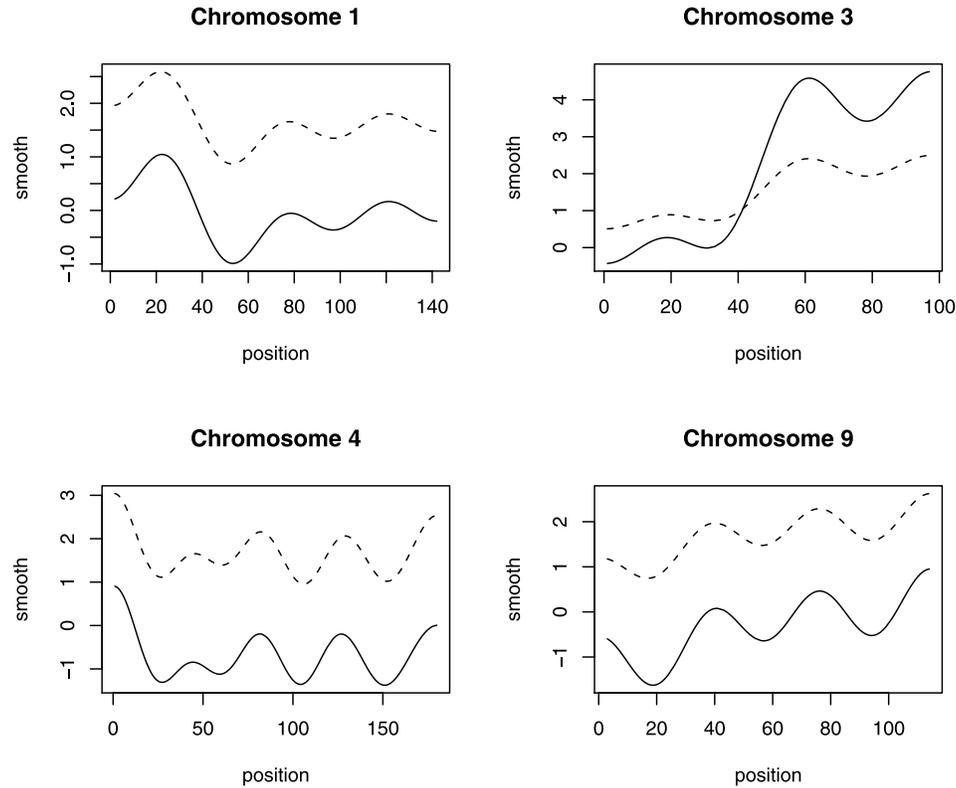

Fig 2. *Fourier series smooths. In each plot, the solid line is a smooth of the raw data, and the dashed line is a smooth of ranks. The smooths are scaled so as to estimate the functions that determine the power of the ordinary and rank-based OS tests.*

The proposed tests can also be applied to test the fit of linear models. In this setting the tests would be applied to residuals, and hence would only be asymptotically distribution free. The bootstrap is another important tool for approximating the sampling distribution of test statistics. An interesting problem for future research is a comparison of bootstrap and rank-based nonparametric lack-of-fit tests. A key question is whether or not large sample rank tests have smaller level error in finite samples than do bootstrap tests.

**Appendix**

Here we provide proofs of Theorems 2.1, 2.2 and 2.3. Throughout the proofs, generic positive (and finite) constants are denoted $C_1, C_2, \ldots$.

*A.1. Proof of Theorem 2.1*

Define $Z_{jn} = \sqrt{24n}\hat{\phi}_j$, $j = 1, \ldots, n-1$, in which case

$$R_n = \max_{1 \leq j < n} \frac{1}{j} \sum_{i=1}^{j} Z_{in}^2.$$

Now define $V_j = F(\epsilon_j)$, $j = 1, \ldots, n$, and note that $V_1, \ldots, V_n$ is a random sample from the $U(0,1)$ distribution. Then $Z_{jn}$ may be expressed

$$Z_{jn} = \sqrt{24n}(\bar{\phi}_j + \delta_j),$$

where

$$\bar{\phi}_j = \frac{1}{n}\sum_{i=1}^{n} V_i \cos(\pi j x_i) \quad \text{and} \quad \delta_j = \frac{1}{n}\sum_{i=1}^{n}(U_i - V_i)\cos(\pi j x_i).$$

Since $F$ is a monotone transformation, $U_i = R(V_i)/(n+1)$, $i = 1, \ldots, n$.

The statistic $R_n$ may be expressed as

$$R_n = \max_{1 \le j < n} \frac{1}{j} \sum_{i=1}^{j} 24n(\bar{\phi}_i^2 + 2\bar{\phi}_i \delta_i + \delta_i^2).$$

At this point we observe that the $U(0,1)$ distribution satisfies the conditions of Theorem 7.2, pp. 168–169, of Hart [3], and therefore

(A.1) $$\max_{1 \le j < n} \frac{1}{j} \sum_{i=1}^{j} 24n\bar{\phi}_i^2 \xrightarrow{\mathcal{D}} T,$$

where $T$ has distribution $G$, as defined in (2.5). Using the following lemma, the proof will then be complete if we can show that

(A.2) $$\max_{1 \le j < n} \frac{1}{j} \sum_{i=1}^{j} n|2\bar{\phi}_i \delta_i + \delta_i^2| \xrightarrow{P} 0.$$

**Lemma A.1.** *For any real numbers $A_1, B_1, \ldots, A_m, B_m$ and any positive number $\epsilon$,*

$$\left| \max_{1 \le j \le m} (A_j + B_j) - \max_{1 \le j \le m} A_j \right| \le \epsilon$$

*whenever $\max_{1 \le j \le m} |B_j| \le \epsilon$.*

*Proof.* Defining $j_1$ and $j_2$ to be such that

$$A_j + B_j \le A_{j_1} + B_{j_1} \quad \text{for } j = 1, \ldots, m,$$

and

$$A_j \le A_{j_2} \quad \text{for } j = 1, \ldots, m,$$

we have

$$\max_{1 \le j \le m}(A_j + B_j) - \max_{1 \le j \le m} A_j = A_{j_1} + B_{j_1} - A_{j_2}.$$

Obviously

$$A_{j_1} + B_{j_2} \le A_{j_2} + B_{j_2} \le A_{j_1} + B_{j_1},$$

and hence

$$B_{j_2} \le A_{j_1} + B_{j_1} - A_{j_2} \le B_{j_1},$$

which proves the result. □

Using Lemma A.1, result (A.1) and Holder's inequality, (A.2) will be proven if we can show that

$$\xi_n \equiv \max_{1 \leq j < n} \frac{1}{j} \sum_{i=1}^{j} n\delta_i^2 \xrightarrow{P} 0.$$

Now

$$\begin{aligned}
P(\xi_n > \epsilon) &= P\left(\bigcup_{j=1}^{n-1}\left\{\frac{1}{j}\sum_{i=1}^{j}\delta_i^2 > \frac{\epsilon}{n}\right\}\right) \\
&\leq \sum_{j=1}^{n-1} P\left(\sum_{i=1}^{j}\delta_i^2 - E_{jn} > \frac{j\epsilon}{n} - E_{jn}\right),
\end{aligned}$$

where $E_{jn} = E\sum_{i=1}^{j}\delta_i^2$. As we shall see shortly, $E_{jn} < j\epsilon/n$ for all $j$ and all $n$ sufficiently large, and so

$$\begin{aligned}
&\sum_{j=1}^{n-1} P\left(\sum_{i=1}^{j}\delta_i^2 - E_{jn} > \frac{j\epsilon}{n} - E_{jn}\right) \\
&\leq \sum_{j=1}^{n-1} P\left(\left[\sum_{i=1}^{j}\delta_i^2 - E_{jn}\right]^2 > \left(\frac{j\epsilon}{n} - E_{jn}\right)^2\right) \\
(A.3) \quad &\leq \sum_{j=1}^{n-1} \operatorname{Var}\left(\sum_{i=1}^{j}\delta_i^2\right)\left(\frac{j\epsilon}{n} - E_{jn}\right)^{-2}.
\end{aligned}$$

Our next step is to determine $E_{jn}$. To this end, define $\Delta_j = V_j - U_j$, $j = 1, \ldots, n$, and note that

$$\begin{aligned}
E_{jn} &= \frac{1}{n^2}\operatorname{Var}(\Delta_1)\sum_{i=1}^{j}\sum_{r=1}^{n}\cos^2(\pi i x_r) \\
(A.4) \quad &+ \frac{1}{n^2}\operatorname{Cov}(\Delta_1, \Delta_2)\sum_{i=1}^{j}\sum_{r=1}^{n}\sum_{s \neq r}\cos(\pi i x_r)\cos(\pi i x_s).
\end{aligned}$$

We now need the following fundamental properties of the cosine basis:

(A.5) $$\sum_{r=1}^{n}\cos(\pi i x_r) = 0, \quad i = 1, \ldots, n-1,$$

and

(A.6) $$\sum_{r=1}^{n}\cos^2(\pi i x_r) = n/2, \quad i = 1, \ldots, n-1.$$

Applying these two properties to (A.4) yields

$$E_{jn} = \frac{j}{2n}\left[\operatorname{Var}(\Delta_1) - \operatorname{Cov}(\Delta_1, \Delta_2)\right].$$

Using basic distributional properties of order statistics and ranks, it is straightforward to show that

$$E_{jn} = \frac{j}{24n^2} + jO(n^{-3}).$$

The bound (A.3) thus satisfies

$$\sum_{j=1}^{n-1} \text{Var}\left(\sum_{i=1}^{j} \delta_i^2\right) \left(\frac{j\epsilon}{n} - E_{jn}\right)^{-2} \leq C_1 n^2 \sum_{j=1}^{n-1} \text{Var}\left(\sum_{i=1}^{j} \delta_i^2\right) \frac{1}{j^2},$$

and the proof will be complete if we can show that

$$\text{Var}\left(\sum_{i=1}^{j} \delta_i^2\right) \leq C_2 \frac{j^2}{n^4}$$

for all $j$ and $n$.

Now,

$$\sum_{i=1}^{j} \delta_i^2 = \frac{1}{n^2} \sum_{r=1}^{n} \sum_{s=1}^{n} \Delta_r \Delta_s \sum_{i=1}^{j} \cos(\pi i x_r) \cos(\pi i x_s)$$

$$= \frac{1}{n^2} \sum_{r=1}^{n} \sum_{s=1}^{n} \Delta_r \Delta_s w_{rs},$$

where $w_{rs} = \sum_{i=1}^{j} \cos(\pi i x_r) \cos(\pi i x_s)$ and for ease of notation we suppress the dependence of $w_{rs}$ on $j$. We have

$$\text{Var}\left(\sum_{i=1}^{j} \delta_i^2\right) = \frac{1}{n^4} \sum_{r=1}^{n} \sum_{s=1}^{n} \sum_{u=1}^{n} \sum_{v=1}^{n} w_{rs} w_{uv} E(\Delta_r \Delta_s \Delta_u \Delta_v) - E_{jn}^2$$

$$\equiv \frac{1}{n^4} S_n - E_{jn}^2,$$

and since we have already established that $E_{jn}^2 \leq C_3 j^2/n^4$, we need only consider $S_n/n^4$. The quantity $S_n$ may be expressed

$$S_n = E(\Delta_1 \Delta_2 \Delta_3 \Delta_4) \sum_r \sum_s \sum_u \sum_v w_{rs} w_{uv}$$

$$+ E(\Delta_1^2 \Delta_2 \Delta_3) \left[ 2 \sum_r \sum_s \sum_u w_{rs} w_{uu} + 4 \sum_r \sum_s \sum_u w_{rs} w_{ru} \right]$$

$$+ E(\Delta_1^2 \Delta_2^2) \left[ \sum_r \sum_s w_{rr} w_{ss} + 2 \sum_r \sum_s w_{rs} w_{rs} \right]$$

(A.7) $\qquad + E(\Delta_1 \Delta_2^3) 4 \sum_r \sum_s w_{rs} w_{ss} + E(\Delta_1^4) \sum_r w_{rr}^2,$

where each sum extends over distinct indices only; for example, $\sum_r \sum_s \sum_u$ is a summation over ordered triples $(r, s, u)$ such that no two of $r, s, u$ are the same.

Each of the sums in the last expression can be shown to be smaller in absolute value than $C_4 j^2 n^2$. Consider, for example, $\sum_r \sum_s \sum_u \sum_v w_{rs} w_{uv}$, the innermost sum of which is

$$\sum_v w_{rs} w_{uv} = w_{rs} \sum_{i=1}^{j} \cos(\pi i x_u) \sum_v \cos(\pi i x_v)$$

$$= -w_{rs} \sum_{i=1}^{j} \cos(\pi i x_u) \left[\cos(\pi i x_r) + \cos(\pi i x_s) + \cos(\pi i x_u)\right],$$

where we have used (A.5). Using in addition (A.6), we have

$$\sum_u \sum_v w_{rs} w_{uv} = 2w_{rs} \sum_{i=1}^{j} \left[\cos^2(\pi i x_r) + \cos^2(\pi i x_s) + \cos(\pi i x_r)\cos(\pi i x_s)\right] - w_{rs}\frac{nj}{2}.$$

Continuing in this way, we finally see that

$$\sum_r \sum_s \sum_u \sum_v w_{rs} w_{uv} = \frac{(nj)^2}{4} + \frac{jn^2}{2} - 6\sum_{r=1}^{n}\left[\sum_{i=1}^{j} \cos^2(\pi i x_r)\right]^2,$$

which is positive and smaller than $C_4 j^2 n^2$ for all $n$ sufficiently large.

We turn now to the expected values in (A.7), each of which is bounded in absolute value by

$$E(\Delta_1^4) = \frac{1}{n}\sum_{i=1}^{n} E\left(U_{(i)} - \frac{i}{n+1}\right)^4$$

$$= \frac{1}{n}\sum_{i=1}^{n} \frac{1}{B(i, n-i+1)} \sum_{j=0}^{4} \binom{4}{j}(-1)^j \left(\frac{i}{n+1}\right)^j \times B(i+4-j, n-i+1),$$

where $B(\cdot, \cdot)$ denotes the beta function. After some straightforward but tedious calculations, it is seen that

$$E(\Delta_1^4) = \frac{1}{10}\cdot\frac{1}{n^2} + O\left(\frac{1}{n^3}\right).$$

Combining this result with previous ones, we have $\mathrm{Var}(\sum_{i=1}^{j} \delta_i^2) \leq C_2 j^2/n^4$, as required, and the proof is complete.

### A.2. Proof of Theorem 2.2

Let $E_{jn}$ denote $E\hat{\phi}_j$, $j = 1, \ldots, n-1$, and define $\delta_j = \hat{\phi}_j - E_{jn}$, $j = 1, \ldots, n-1$. Then

$$\hat{\phi}_j = \delta_j + \frac{h(0)\beta_j}{\sqrt{n}} + \left[E_{jn} - \frac{h(0)\beta_j}{\sqrt{n}}\right],$$

and it suffices to show that

$$\max_{1 \leq m < n} \frac{1}{m}\sum_{j=1}^{m} \left(\sqrt{24n}\delta_j + \sqrt{24}h(0)\beta_j\right)^2 \to$$

(A.8)
$$\sup_{m \geq 1} \frac{1}{m}\sum_{j=1}^{m}\left(Z_j + \sqrt{24}h(0)\beta_j\right)^2$$

in distribution as $n \to \infty$, and that

(A.9)
$$\lim_{n\to\infty} \max_{1\leq m<n} \frac{1}{m}\sum_{j=1}^{m}(\sqrt{n}E_{jn} - h(0)\beta_j)^2 = 0.$$

The proof of (A.8) is virtually the same as the proof of Theorem 2.1 and hence omitted.

To prove (A.9), first notice that

$$\max_{1 \leq m < n} \frac{1}{m} \sum_{j=1}^{m} (\sqrt{n} E_{jn} - h(0)\beta_j)^2 \leq$$

$$\max \left[ \max_{1 \leq m < \sqrt{n}} \frac{1}{m} \sum_{j=1}^{m} (\sqrt{n} E_{jn} - h(0)\beta_j)^2, \frac{1}{\sqrt{n}} \sum_{j=1}^{n-1} (\sqrt{n} E_{jn} - h(0)\beta_j)^2 \right].$$

We have

$$(A.10) \quad \sum_{j=1}^{n-1} (\sqrt{n} E_{jn} - h(0)\beta_j)^2 \leq \left[ \left( \sum_{j=1}^{n-1} n E_{jn}^2 \right)^{1/2} + h(0) \left( \sum_{j=1}^{n-1} \beta_j^2 \right)^{1/2} \right]^2.$$

Defining

$$\mu_n(x) = \frac{1}{n} \sum_{i=1}^{n} H\left( \frac{\beta(x) - \beta(x_i)}{\sqrt{n}} \right) \quad \forall\, x$$

and $\bar{\mu}_n = \sum_{i=1}^{n} \mu_n(x_i)/n$, it follows from orthogonality of the cosines that

$$(A.11) \quad \sum_{j=1}^{n-1} E_{jn}^2 = \frac{1}{2n} \sum_{i=1}^{n} (\mu_n(x_i) - \bar{\mu}_n)^2.$$

Now, $\mu_n(x_i) = H(0) + n^{-3/2} \sum_{k=1}^{n} (\beta(x_i) - \beta(x_k)) h(\eta_{ikn})$, where $\eta_{ikn}$ is a number between 0 and $(\beta(x_i) - \beta(x_k))/\sqrt{n}$. It follows that

$$\sqrt{n}(\mu_n(x_i) - \bar{\mu}_n) = \frac{1}{n} \sum_{k=1}^{n} (\beta(x_i) - \beta(x_k)) h(\eta_{ikn})$$

$$(A.12) \qquad\qquad - \frac{1}{n^2} \sum_{j=1}^{n} \sum_{k=1}^{n} (\beta(x_j) - \beta(x_k)) h(\eta_{jkn}).$$

Since $\beta$ and $h$ are both bounded functions, (A.11) and (A.12) imply that $n \sum_{j=1}^{n-1} E_{jn}^2$ can be bounded by the same constant for all $n$. By the piecewise continuity of $\beta$ on $[0, 1]$, $\sum_{j=1}^{n-1} \beta_j^2$ has a finite limit, which, using (A.10), finishes the proof that $\sum_{j=1}^{n-1} (\sqrt{n} E_{jn} - h(0)\beta_j)^2 / \sqrt{n}$ tends to 0 as $n \to \infty$.

Finally, we need to show that

$$(A.13) \quad \lim_{n \to \infty} \max_{1 \leq m < \sqrt{n}} \frac{1}{m} \sum_{j=1}^{m} (\sqrt{n} E_{jn} - h(0)\beta_j)^2 = 0.$$

We may write $\sqrt{n} E_{jn} - h(0)\beta_j = A_{nj} + B_{nj} + C_{nj}$, where

$$A_{nj} = h(0) \left[ \frac{1}{n} \sum_{i=1}^{n} \beta(x_i) \cos(\pi j x_i) - \beta_j \right],$$

$$B_{nj} = \frac{1}{n} \sum_{i=1}^{n} \beta(x_i) \cos(\pi j x_i) \frac{1}{n} \sum_{k=1}^{n} [h(\eta_{ikn}) - h(0)] \quad \text{and}$$

$$C_{nj} = \frac{1}{n^2} \sum_{i=1}^{n} \cos(\pi j x_i) \sum_{k=1}^{n} \beta(x_k) [h(0) - h(\eta_{ikn})].$$

Using the boundedness of $\beta$ and Lipschitz continuity of $h$ near 0, we have, for all $n$ sufficiently large,
$$|h(\eta_{ikn}) - h(0)| \leq \frac{C_1}{\sqrt{n}} \quad \text{for all } i \text{ and } k.$$

This fact along with boundedness of $\beta$ and the cosines implies that $|B_{nj} + C_{nj}| \leq C_2/\sqrt{n}$ for all $n$ sufficiently large and each $j$. Now, using the second of the two conditions in Theorem 2.2, it it straightforward to show that $|A_{nj}| \leq C_3 j/n$ for all $j$ and $n$.

Combining the results in the previous paragraph, it is clear that
$$\max_{1 \leq m < \sqrt{n}} \frac{1}{m} \sum_{j=1}^{m} (\sqrt{n} E_{jn} - h(0)\beta_j)^2 = O\left(\frac{1}{n}\right)$$

as $n \to \infty$, which completes the proof of Theorem 2.2.

### *A.3. Proof of Theorem 2.3*

The proof of Theorem 2.3 is simple given result (2.7) and the proof of Theorem 2.2. Let the alternative function be $\beta(x)/\sqrt{n}$, where $n$ is the sample size of the ordinary OS test, and let the sample size, $n^*$, for the rank-based OS test be
$$n^* = \left[\frac{n}{12\sigma^2 h^2(0)}\right],$$

where $[z]$ = greatest integer less than or equal to $z$. If we can show that the limiting power of the rank-based OS test with sample size $n^*$ is equal to (2.7), then the proof is complete.

We have
$$24 n^* \hat{\phi}_j^2 = \left[\sqrt{24 n^*}\delta_j + \frac{\sqrt{2}\beta_j}{\sigma} + \sqrt{24 n^*}\left(E_{jn^*} - \frac{h(0)\beta_j}{\sqrt{n}}\right) \right.$$
$$\left. + \frac{\sqrt{24 n^*} h(0) \beta_j}{\sqrt{n}} - \frac{\sqrt{2}\beta_j}{\sigma}\right]^2,$$

where
$$\hat{\phi}_j = \frac{1}{n^*} \sum_{i=1}^{n^*} U_i \cos(\pi j x_i), \quad E_{jn^*} = E\hat{\phi}_j \quad \text{and} \quad \delta_j = \hat{\phi}_j - E_{jn^*}.$$

The remainder of the proof is virtually the same as that of Theorem 2.2. A key point is that, to first order, $E\hat{\phi}_j$ is still $h(0)\beta_j/\sqrt{n}$, which is true because the factor $\sqrt{n}$ in $h(0)\beta_j/\sqrt{n}$ derives from the local alternative and not the sample size $n^*$.